\documentclass[12pt,reqno]{amsart}
\usepackage{amsmath}
\usepackage{amsxtra}
\usepackage{amscd}
\usepackage{amsthm}
\usepackage{amsfonts}
\usepackage{amssymb}
\usepackage{eucal}

\newtheorem{theorem}{Theorem}[section]

\newtheorem{cor}[theorem]{Corollary}
\newtheorem{lem}[theorem]{Lemma}

\theoremstyle{definition}

\theoremstyle{remark}

\theoremstyle{remark}

\numberwithin{equation}{section}

\newcommand{\nc}{\newcommand}
\nc{\on}{\operatorname}
\nc{\ch}{\mbox{ch}}
\nc{\Z}{{\mathbb Z}}
\nc{\C}{{\mathbb C}}
\nc{\R}{{\mathbb R}}
\nc{\pone}{{\mathbb C}{\mathbb P}^1}
\nc{\pa}{\partial}
\nc{\F}{{\mathcal F}}
\nc{\arr}{\rightarrow}
\nc{\larr}{\longrightarrow}
\nc{\al}{\alpha}
\nc{\ri}{\rangle}
\nc{\lef}{\langle}
\nc{\W}{{\mathcal W}}
\nc{\la}{\lambda}
\nc{\ep}{\epsilon}

\nc{\Om}{\Omega}
\nc{\su}{\widehat{{\mathfrak sl}}_2}
\nc{\sw}{{\mathfrak s}{\mathfrak l}}

\nc{\g}{{\mathfrak g}}
\nc{\h}{{\mathfrak h}}
\nc{\n}{{\mathfrak n}}
\nc{\N}{\widehat{\n}}
\nc{\G}{\widehat{\g}}
\nc{\De}{\Delta_+}
\nc{\gt}{\widetilde{\g}}
\nc{\Ga}{\Gamma}
\nc{\one}{{\mathbf 1}}
\nc{\z}{{\mathfrak Z}}
\nc{\zz}{{\mathcal Z}}
\nc{\Hh}{{\mathcal H}_\beta}
\nc{\qp}{q^{\frac{k}{2}}}
\nc{\qm}{q^{-\frac{k}{2}}}
\nc{\La}{\Lambda}
\nc{\wt}{\widetilde}
\nc{\qn}{\frac{[m]_q^2}{[2m]_q}}
\nc{\cri}{_{\on{cr}}}
\nc{\kk}{h^\vee}
\nc{\sun}{\widehat{\sw}_N}
\nc{\hh}{\widehat{\mathfrak h}}
\nc{\HH}{{\mathcal H}_{q,t}}
\nc{\ca}{\wt{{\mathcal A}}_{h,k}(\sw_2)}
\nc{\gl}{\widehat{{\mathfrak g}{\mathfrak l}}_2}
\nc{\el}{\ell}
\nc{\s}{{\mathbf s}}
\nc{\bi}{\bibitem}
\nc{\om}{\omega}
\nc{\WW}{\W_\beta}
\nc{\scr}{{\mathbf S}}
\nc{\ab}{{\mathbf a}}
\nc{\rr}{r}
\nc{\ol}{\overline}
\nc{\con}{qt^{-1} + q^{-1}t}
\nc{\den}{q^{\el-1} t^{-\el+1}+ q^{-\el+1} t^{\el-1}}
\nc{\ds}{\displaystyle}
\nc{\B}{B}
\nc{\A}{{\mathbb A}}
\nc{\GG}{{\mathcal G}}
\nc{\UU}{{\mathcal U}}
\nc{\MM}{{\mathcal M}}
\nc{\CC}{{\mathcal C}}
\nc{\GL}{{}^L G}
\nc{\dzz}{\frac{dz}{z}}
\nc{\Res}{\on{Res}}
\nc{\rep}{{\mathcal R}ep \;}
\nc{\uqg}{U_q \G}
\nc{\uqgg}{U_q \g}
\nc{\Fq}{{\mathbb F}_q}

\nc{\stimes}{\ltimes}
\nc{\K}{\hat{\mathcal K}}
\nc{\Ql}{\ol{\mathbb Q}_\ell}

\nc{\ga}{\gamma}
\nc{\PL}{{}^L P}
\nc{\E}{\mc E}
\nc{\mc}{\mathcal}
\nc{\mbf}{\mathbf}
\nc{\bb}{{\mathfrak b}}
\nc{\OO}{{\mc O}}
\nc{\Po}{{\mc P}}
\nc{\V}{{\mc V}}
\nc{\yy}{{\mc Y}}
\nc{\M}{\mathcal M}
\nc{\Coh}{{{\mathcal C}oh}}
\nc{\Cohn}{\Coh_n}
\nc{\f}{{\mathcal F}}
\nc{\si}{_E}
\nc{\Gaf}{{\mathbb G}_{a,\Fq}}
\nc{\KK}{{\mathfrak k}}

\nc{\PCr}{{ \bs P  (\C[x])^r   }}
\nc{\PCN}{{ \bs P  (\C[x])^N   }}
 
\nc{\sN}{sl_{2N+1}}

\nc{\Pzr}{{ \bs P(\C((x-z)))^r}}

\nc{\PzN}{{ \bs P(\C((x-z)))^N}}

\newcommand{\di}{ \frac d{dx}}

\newcommand{\bean}{\begin{eqnarray}}
\newcommand{\eean}{\end{eqnarray}}
\newcommand{\be}{\begin{displaymath}}
\newcommand{\ee}{\end{displaymath}}
\newcommand{\bea}{\begin{eqnarray*}}   
\newcommand{\eea}{\end{eqnarray*}}
\newcommand{\bs}{\boldsymbol}
\newcommand{\Ref}[1]{{$($\ref{#1}$)$}}

\begin{document}

\title{Differential equations for Jacobi-Pi\~neiro polynomials}
\author[E. Mukhin and A. Varchenko]
{E. Mukhin and A. Varchenko}
\thanks{
Research of A.V. is supported in part by NSF grant DMS-0244579.}
\address{E.M.: Department of Mathematical Sciences, Indiana University -
Purdue University Indianapolis, 402 North Blackford St, Indianapolis,
IN 46202-3216, USA, \newline mukhin@math.iupui.edu}
\address{A.V.: Department of Mathematics, University of North Carolina 
at Chapel Hill, Chapel Hill, NC 27599-3250, USA, anv@email.unc.edu}
\maketitle

\begin{abstract}
For $r\in \Z_{\geq 0}$, we present  a Fuchsian linear
differential operator of order $r+1$
with three singular points at $0, 1, \infty$.
This operator annihilates the $r$-multiple Jacobi-Pi\~neiro
polynomial.
\end{abstract}

\section{Introduction}
Let $r$ be a natural number.
Consider a Fuchsian differential operator
$$
D\ = \  \sum_{i=0}^{r+1} \,c_i(x)\, \frac {d^i}{dx^i}\ 
$$
with singular points at $z_1,\dots,z_n, \infty$ and with kernel
consisting of polynomials only. An interest to such operators had
arisen recently in relation with the Bethe ansatz method in the Gaudin
model, where such operators were used to construct eigenvectors of the
Gaudin Hamiltonians, see \cite{ScV}, \cite{MV1}-\cite{MV3},
\cite{MTV1}, \cite{MTV2}.

In the Gaudin
model, one considers the tensor product $M = M_1\otimes \dots \otimes M_n$
of finite dimensional irreducible $\frak{gl}_{r+1}$-modules,
located respectively at $z_1,\dots,z_n$.  The module $M_s$,
sitting at $z_s$, is determined by the exponents of $D$ at $z_s$.  One
constructs $r+1$ one-parameter families of commuting linear operators
$H_i(x) : M \to M$, $i=1,\dots,r+1$, acting on $M$ and called the
Gaudin Hamiltonians. The problem is to construct eigenvectors and eigenvalues
of the Gaudin Hamiltonians. 

It turns out, that having the kernel of the differential operator
$D$, i.e.  the $r+1$-dimensional vector space of polynomials, 
one constructs (under certain conditions) an eigenvector $v_D \in M$ of
the Gaudin Hamiltonians with corresponding eigenvalues being
the coefficients of $D$,
$$
H_i(x)\, v_D\ =\ c_i(x)\,v_D\ , 
\qquad 
 i = 1,\dots, r+1 \ .
$$

The Bethe ansatz idea is to construct all eigenvectors of the Gaudin Hamiltonians
by choosing different operators $D$ with the same singular points and exponents.

This philosophy motivates the detailed study of 
Fuchsian operators with prescribed singular points,
exponents, and polynomial kernels.

The important model case is the study of operators with three singular points
$0,1,\infty$.

The operators with special exponents $0,k+1,k+2, \dots,k+r$
at $x=1$ and arbitrary exponents at $x=0, \infty$ were studied in
\cite{MV2}. It was discovered in \cite{MV2} that the kernel of such a
differential operator  consists of
Jacobi-Pi\~neiro polynomials, a special type of multiple orthogonal
polynomials, see  Lemma 4.4 in \cite{MV2}.

This appearance of orthogonal polynomials in the Bethe ansatz
constructions helped us in \cite{MV2} study eigenvectors of the Gaudin
Hamiltonians.

In this short paper, we give an example of a reverse implication, namely, 
that the Bethe ansatz considerations may be useful in
studying  orthogonal polynomials.  We construct a Fuchsian
differential operator with singular points at $0,1,\infty$ 
annihilating the Jacobi-Pi\~neiro polynomial, see the precise statement 
and the discussion of the result in Section 5. Such an operator can be used in 
studying the Jacobi-Pi\~neiro polynomials.

\medskip

We thank referees for helping to improve the exposition.

\section{Jacobi-Pi\~neiro polynomials}

Let $l_1,\dots,l_r$ be integers such that $l_1\geq \dots \geq l_r\geq 0$.
Let $m_1,\dots,m_r$ and $k$ be negative real numbers. We use the notation 
$\bs m=(m_1,\dots,m_r)$, $\bs l=(l_1,\dots,l_r)$.

The Jacobi-Pi\~neiro polynomial \cite{P} is the 
unique monic polynomial of degree $l_1$ 
whose coefficients are rational functions of $\bs m,\bs l,k$ and 
which is orthogonal to functions
\bea
\underbrace{1,x,\dots,x^{l_1-l_2-1}}_{l_1-l_2}\ ,
\underbrace{x^{-m_2-1},x^{-m_2},\dots, x^{-m_2+l_2-l_3-2}}_{l_2-l_3}\ ,
\dots, \underbrace{x^{-\sum_{i=2}^{r}m_i-r+1},\dots,
x^{-\sum_{i=2}^{r}m_i-r+l_r}}_{l_r}
\eea
with respect to the scalar product given by
\be
(f(x),g(x))=\int_0^1f(x)g(x)(x-1)^{-k-1}x^{-m_1-1}dx.
\ee
We denote the Jacobi-Pi\~neiro polynomial  by $P_{\bs m,\bs l,k}(x)$.

If $l_2=l_3=\dots=l_r=0$, then the Jacobi-Pi\~neiro polynomial is
the classical Jacobi polynomial $P_l^{(\al,\beta)}(x)$ on
interval $[0,1]$ with $l=l_1$, $\al=-k-1$, $\beta=-m_1-1$.

The Jacobi-Pi\~neiro polynomial may be given by the Rodrigues-type formula, 
see \cite{ABV}:
\begin{align}\label{Rodr}
P(& \bs m,\bs l,k) \ =\ c\, (x-1)^{k+1}\,x^{{}\,\sum_{i=1}^rm_i\,-\,r}\ \times
\\
& \times  \frac{d^{l_r-l_{r+1}}}{d x^{l_r-l_{r+1}}}\
   x^{l_r-l_{r+1}-m_r-1}\frac{d^{l_{r-1}-l_r}}{d
     x^{l_{r+1}-l_r}}\ \dots\
   x^{l_2-l_3-m_2-1}\frac{d^{l_1-l_2}}{d x^{l_{1}-l_2}}\
   \left( x^{l_1-l_2-m_1-1}(x-1)^{l_1-k-1}\right), \notag
 \end{align}
where $c$ is a nonzero constant.

The coefficients of the Jacobi-Pi\~neiro polynomial $P_{\bs m,\bs
l,k}(x)$ are rational functions of $\bs m,\bs l,k$ and therefore the
polynomial $P_{\bs m,\bs l,k}(x)$ is well defined for almost all
complex $m_1,\dots,m_r,k$.

\section{Spaces of polynomials the first and second type } 
We describe remarkable spaces of polynomials 
which contain Jacobi-Pi\~neiro polynomials. 
See \cite{MV2} for the relation of these spaces to the Bethe Ansatz method.

Parameters $(\bs m,\bs l,k)$ are called {\it consistent} if all 
$m_i,l_i$ and $k$ are nonnegative integers satisfying
\bean\label{consistent}
k\geq l_1\geq l_2\geq \dots\geq l_r\geq 0, \qquad  l_{s}-l_{s+1}\leq m_{s}
\qquad 
(s=1,\dots,r)\ .
\eean

Let  $(\bs m,\bs l,k)$ be consistent. We use the convention:
\be 
l_0=k, \qquad l_{r+1}=0\ .
\ee

We call a complex $r+1$-dimensional
vector space of polynomials $V(\bs m,\bs l,k)\subset \C[x]$ {\it
the space of polynomials of the first type associated to $(\bs m,\bs
l,k)$} if the space satisfies the following two conditions:
\begin{itemize}
\item  The space  $V(\bs m,\bs l,k)$ has a basis of the form
\bean
\label{b1}
\{v_0(\bs m, \bs l,k),\, v_1(\bs m, \bs l,k)\,x^{m_1+1},\,
v_2(\bs m, \bs l,k)\,x^{m_1+m_2+2},\,
\dots,\, v_r(\bs m, \bs l,k)\,x^{\sum_{i=1}^rm_i\,+\,r}\},
\eean
where for $i=0,\dots, r$, the polynomial $v_i(\bs m, \bs l,k)\in\C[x]$ is a 
monic polynomial of degree $k-l_i+l_{i+1}$\,.
\item 
If a polynomial $p \in V(\bs m,\bs l,k)$ vanishes at 
$x=1$, then the multiplicity of zero at $x=1$ is at  least $k+1$. 
\end{itemize}
It is easy to see that the basis polynomials in \Ref{b1}
have increasing degrees.

Below we will show that for any consistent parameters $(\bs m,\bs l,k)$, 
there exists a unique space of polynomials of the first type associated to
$(\bs m,\bs l,k)$, see Theorem \ref{unique}. 
Moreover, we will show that this space contains the 
Jacobi-Pi\~neiro polynomial $P(\bs m,\bs l,k)$, see Lemma \ref{P in V}.

We call a complex $r+1$-dimensional
vector space of polynomials $U(\bs m,\bs l,k) \subset \C[x]$ 
\linebreak
{\it
the space of polynomials of the second type associated to $(\bs m,\bs
l,k)$} if the space satisfies the following two conditions:
\begin{itemize}
\item 
The space  $U(\bs m,\bs l,k)$ has a basis of the form
\bean
\label{b2}
\{u_0(\bs m, \bs l,k),\, u_1(\bs m, \bs l,k)\,x^{m_r+1}\!,\,
u_2(\bs m, \bs l,k)\,x^{m_r+m_{r-1}+2}\!,
\dots, 
u_r(\bs m, \bs l,k)\,x^{\sum_{i=1}^rm_i\,+\,r}\},
\eean
where for $i=0,\dots,r$, the polynomial $u_i(\bs m, \bs l,k)\in\C[x]$ is a monic
polynomial of degree $l_{r-i}-l_{r-i+1}$\,.
\item 
There exists a nonzero polynomial $p\in U(\bs m,\bs l,k)$ which has 
zero at $x=1$ of order $k+r$.
\end{itemize}
It is easy to see that the basis polynomials in \Ref{b2}
have increasing degrees.

Below we will show that for any consistent parameters $(\bs m,\bs l,k)$,
there exists a unique space of polynomials of the second type associated to
$(\bs m,\bs l,k)$, see Theorem \ref{unique}.

The spaces of the first type and of the second type are dual in the 
sense of \cite{MV1} which we now describe.

Define an $r$-tuple $\bs T=(T_1,\dots,T_r)$ of polynomials in $x$ by 
\bean\label{T}
T_1\ =\ (x-1)^k\,x^{m_1}\ , 
\qquad 
T_i\ =\ x^{m_i} 
\qquad (\,i = 2 , \dots , r\,)\ .
\eean
For functions $f_1,\dots,f_s$ of $x$,  
{\it the Wronskian  $W(f_1,\dots,f_s)$} is defined by
\be
W(f_1, \dots , f_s)=\det\left(\frac{d^{i}}{d x^{i}}\ 
f_{j}\right)_{i,j=1,\dots,s}.
\ee
For functions $f_1, \dots, f_s$ of $x$,
define {\it the divided Wronskians} 
$W^\dagger_V(f_1, \dots , f_s)$ and 
\linebreak
$W^\dagger_U(f_1, \dots , f_s)$ by
\bea
&&
W^\dagger_V(f_1, \dots , f_s)\ 
=
\ W(f_1,\dots,f_s)\ T_1^{1-s} \,T_2^{2-s} \dots T_{s-1}^{-1}\ ,
\\
&&
W^\dagger_U(f_1, \dots , f_s)\ =\ W(f_1, \dots , f_s)\ 
T_r^{1-s} \,T_{r-1}^{2-s} \dots T_{r-s+2}^{-1}\ .
\eea

\begin{lem}\label{dual}
Let  $(\bs m,\bs l,k)$ be consistent parameters.

Let $V$ be a space of the first type associated to $(\bs m,\bs l,k)$.
Then the space 
\bea
U\ =\ \{\,W^\dagger_V(f_1,\dots,f_r)\,, \ {} \,f_1, \dots , 
f_r \in V \}\ 
\eea
is a space of polynomials of the second type associated to  $(\bs m,\bs l,k)$.

Let $U$ be a space of the second type associated to $(\bs m,\bs l,k)$.
Then the space 
\bea
V\ =\ \{\,W^\dagger_U(f_1,\dots,f_r)\,, \ {}\, f_1, \dots , 
f_r \in U\, \}\ 
\eea
is a space of polynomials of the first type associated to  $(\bs m,\bs l,k)$.
\end{lem}
\begin{proof}
  The lemma follows from the definitions.
\end{proof}

\section{Recursion for spaces $V(\bs m,\bs l,k)$} 
We show the existence of spaces $V(\bs m,\bs l,k)$ of the first type
by constructing them recursively as follows.

Let $m_1,\dots,m_r$ be nonnegative numbers. Let $\bs 0 = (0, \dots , 0)$.
Then clearly the parameters $(\bs m, \bs 0, k=0)$ are consistent.

Introduce the numbers 
$$
e_i\ =\ i\ +\ \sum_{j=1}^i\, m_j\ , 
\qquad
(\, i = 0, \dots, r\, )\ .
$$ 
In particular,  $e_0 = 0$.

\begin{lem}\label{start}
The space 
\be
V(\bs m,\bs 0,0)\ =\ {\rm span}\, \langle\,1=x^{e_0}\,,\
x^{e_1}\,,\ \dots ,\ x^{e_r}\,\rangle
\ee
is a space of the first type associated to $(\bs m,\bs 0,0)$.
\end{lem}
\begin{proof}
The lemma is proved by direct verification.
\end{proof}

For $i=0,1,\dots,r$, introduce the first order linear
differential operators
\begin{align}\label{op1}
D_i(\bs m,\bs l,k)\ =\ 
x(x-1)\,\frac d {dx}\ -\ (k\,+\,\sum_{s=1}^{i}m_s-l_{i}+l_{i+1}+i)(x-1) - k - 1\ .
\end{align}
For $i=1,\dots, r$, let
$\bs 1_i\,=\,(1,\dots,1,1,0,\dots,0)$ be the $r$-tuple 
where we have $i$ ones and $r-i$ zeros.
Let $\bs 1_0\,=\,\bs 0 = (0,\dots,0)$.

For all $i,j\in\{0,1,\dots,r\}$, we have
\be
D_j(\bs m,\bs l+\bs 1_i,k+1)\ 
D_i(\bs m,\bs l,k)\ =\ 
D_i(\bs m, \bs l+\bs 1_j,k+1)\ D_j(\bs m,\bs l,k)\ .
\ee

\begin{lem}\label{V lem}
  Suppose $(\bs m, \bs l, k)$ and $(\bs m, \bs l+\bs 1_i, k+1)$ are
  consistent parameters. Let $V(\bs m,\bs l,k)$ be a space of the
  first type associated to $(\bs m, \bs l, k)$.  Then the space
\be  
 V(\bs m, \bs l+\bs 1_i, k+1)\
=\ \{\,D_i(\bs m,\bs l,k)\,v\ , \ {}\, v\in V(\bs m,\bs l,k)\,\} 
\ee
 is a space of the first type associated to $(\bs m, \bs l+\bs 1_i, k+1)$.
\end{lem}
\begin{proof}
The proof is straightforward.
\end{proof}

\begin{theorem}\label{exist}
Let  $(\bs m, \bs l, k)$ be consistent parameters. Then there exist 
a space $V(\bs m,\bs l,k)$ of the first type and a space $U(\bs m,\bs l,k)$ of the 
second type associated to $(\bs m, \bs l, k)$.
\end{theorem}

\begin{proof}
Let $i(1),\,\dots\,,\,i(k)\,\in \,\{\,0,\dots,r\,\}$\ 
be any sequence of indices such that 
$0$ occurs in the sequence exactly
$k-l_1$ times and every number \ $i\,=\,1, \dots , r$ 
occurs in the sequence exactly $l_i-l_{i+1}$ times.
Then $\bs l\,=\, \sum_{s=1}^k \bs 1_{i(s)}$ and
for every for $j=0,\dots,k$, the tuple
 $(\bs m, \sum_{s=1}^j\bs 1_{i(s)},j)$ forms a 
consistent set of parameters.

Introduce the linear differential operator
\begin{align}\label{DD}
D_{\bs m,\bs l,k}\ =\
  D_{j(k)}(\bs m,\sum_{s=1}^{k-1}\bs 1_{j(s)},k-1)\ \dots\ 
D_{j(2)}(\bs m, \bs 1_{j(1)},1)\ D_{j(1)}(\bs m,\bs 0,0)\ 
\end{align}
of order $k$.

By Lemma \ref{op1}, a space $V(\bs m,\bs l,k)$    
of the first type associated to $(\bs m,\bs l,k)$ 
can be constructed by application of the operator $D_{\bs m,\bs l,k}$
to the space $V(\bs m,\bs 0,0)$ of Lemma \ref{start}. 

A space $U(\bs m,\bs l,k)$  
of the second type associated to $(\bs m,\bs l,k)$
can be constructed from the space of the first type
by the construction of  Lemma \ref{dual}.
\end{proof}

Let  $(\bs m, \bs l, k)$ be consistent parameters. 
Let  $V(\bs m,\bs l,k)$ be the space of the first type associated to
$(\bs m, \bs l, k)$.
Let $v_0(\bs m,\bs l,k)\in V(\bs m,\bs l,k)$ be the monic polynomial of
degree $l_1$. Such a polynomial in $V(\bs m,\bs l,k)$
is unique according to the definition of the space
of the first type associated to $(\bs m,\bs l,k)$.

\begin{lem}
\label{P in V}
The polynomial $v_0(\bs m,\bs l,k)$ is the
  Jacobi-Pi\~neiro polynomial $P(\bs m,\bs l,k)$.
\end{lem}
\begin{proof}
  The polynomial $v_0(\bs m,\bs l,k)$ is obtained by  application
  of the operator $D_{\bs m,\bs l,k}$ to the function $1$. It is a straightforward 
  calculation to check that this formula for $v_0(\bs m,\bs
  l,k)$ coincides with the Rodrigues-type formula for the
  Jacobi-Pi\~neiro polynomial $P(\bs m,\bs l,k)$ in formula \Ref{Rodr}.
\end{proof}

\section{The differential operator $D^\vee_{\bs m,\bs l,k}$.}
Let $(\bs m,\bs l,k)$ be consistent parameters. 
We recall our convention $l_0=k$, and $l_{r+1}=0$.
Consider two sets of  numbers,
\bea
d_i(\bs m,\bs l,k)\ =\
\sum_{s=r+1-i}^r\ m_s-l_{r-i+1}+l_{r-i}+i \ ,
\qquad
a_i(\bs m,\bs l,k)\ =\ \sum_{s=r+1-i}^rm_s+i \ ,
\eea
where  $i\,=\,0 , \dots  , r$, \
and two polynomials in $\al$ whose roots are those numbers,
\bea
d(\al;\bs m,\bs l,k)\ =\ \prod_{i=0}^r(\al-d_i(\bs m,\bs l,k))\ ,
 \qquad
a(\al;\bs m,\bs l,k)\ = \ \prod_{i=0}^r(\al-a_i(\bs m,\bs l,k))\ .
\eea
Define the numbers $A_i$ and $B_i$ for $i=0,\dots,r+1$, as the coefficients 
of the following decompositions
\bean\label{infty}
d(\al;\bs m,\bs l,k)\ =\ A_0(\bs m,\bs l,k)+ \sum_{i=1}^{r+1}\ A_i (\bs m,\bs l,k)\,
\al(\al-1)\dots(\al-i+1)\ ,
\\
a(\al;\bs m,\bs l,k)\ =\ B_0(\bs m,\bs l,k)+ \sum_{i=1}^{r+1}\ B_i (\bs m,\bs l,k)\,
\al(\al-1)\dots(\al-i+1)\ .
\label{zero}
\eean
Clearly we have 
$A_{r+1}(\bs m,\bs l,k) = B_{r+1}(\bs m,\bs l,k) = 1,$\ 
$A_0(\bs m,\bs l,k) = d(0;\bs m,\bs l,k)$, \
$B_0(\bs m,\bs l,k)  = a(0;\bs m,\bs l,k) = 0$.

\medskip
\noindent
{\bf Remark.}
 If $f$ is a polynomial in $\al$ and
$\Delta f(\al) = f(\al+1)-f(\al)$, then
\bea 
 f(\al)\ =\ \sum_{i}\
\frac{\Delta^if(0)}{i!}\,\al(\al-1)\dots(\al-i+1)\ . 
 \eea

\medskip

Introduce the monic linear differential operator
\be 
D^\vee_{\bs m,\bs l,k}\ =\  \sum_{i=0}^{r+1}\ 
\frac {A_i(\bs m,\bs l,k)\,x\, -\, B_i(\bs m,\bs l,k)}
{x^{r+1-i}\,(x-1)} \ \frac {d^i}{dx^i}\ 
\ee
of order $r+1$.

\begin{lem}\label{DveeU}
Let $U(\bs m,\bs l,k)$ be a space of the second type associated to  
$(\bs m,\bs l,k)$. Then $U(\bs m,\bs l,k)$ is the
kernel of $D^\vee_{\bs m,\bs l,k}$\, .
\end{lem}
\begin{proof}
  Let 
$$
\tilde D\ =\ \sum_{i=0}^{r+1} \,c_i(x)\, \frac {d^i}{dx^i}
$$ 
be the monic
  differential operator of order $r+1$ with kernel $U(\bs m,\bs l,k)$.
The operator
  $\tilde D$ is a Fuchsian differential
operator with singular points at $0,1,\infty$. 
The coefficients $c_i$ are rational functions in $x$.
A coefficient $c_i$ may have poles only at $x=0$ and $x=1$ of orders
 at most $i-r-1$ and the degree of $c_i$ at infinity is at most $i-r-1$. 
It is easy to see that
the poles of coefficients $c_i$ at $x=1$ are at most
 simple, cf. for example, formula (5.1) in \cite{MV1}.  
Therefore, the coefficients $c_i$
 can be written in the form 
 \be
 c_i\  =\ \frac {\tilde A_i\,x\ -\ \tilde  B_i}{x^{r+1-i}\,(x-1)}\ .  
 \ee 
 From the characteristic equation for exponents of $\tilde D$ at
 $x=\infty$ and formula \Ref{infty}, we conclude that $\tilde A_i=A_i(\bs m,\bs
 l,k)$.  From the characteristic equation for exponents
of $\tilde D$ at $x=0$ and formula
\Ref{zero}, we
 conclude that $\tilde B_i=B_i(\bs m,\bs l,k)$.
\end{proof}

\begin{theorem}
\label{unique}
The space $V(\bs m,\bs l,k)$ of the first type and the space $U(\bs m,\bs l,k)$ of the 
second type associated to $(\bs m,\bs l,k)$ are unique.
\end{theorem}
\begin{proof}
The space of the second type is unique by Lemma \ref{DveeU}.
The space of the first type is unique by Lemma \ref{dual}, since formulas of
Lemma \ref{dual} allow us recover uniquely the space of the first type from 
the unique space of the second type.

\end{proof}

\medskip

\noindent
{\bf Remark.}
Set 
\bean
\label{REMAR}
T(x)\ =\ (x-1)^{k}\, x^{{} \sum_{i=1}^r i \,m_i }\ .
\eean
It is easy to see that if $f_1,\dots, f_{r+1}$ is a basis of
$U(\bs m,\bs l,k)$, then the Wronskian of 
$f_1,\dots, f_{r+1}$ is equal to $T$ up to multiplication by a number.

\section{Differential equation for the Jacobi-Pi\~neiro polynomial}
In this section we present a linear differential operator of order
$r+1$ annihilating the Jacobi-Pi\~neiro polynomial $P_{\bs m,\bs l,k}(x)$.

Set 
\be
\tau(x)\ =\ (x-1)^k\,x^{\sum_{i=1}^r\,m_i}\ .
\ee

Define  the linear differential operator,
\bea
 D_{\bs m,\bs l,k}\ =\ \tau(x) 
\sum_{i=0}^{r+1}\ (-1)^{r+1+i} \frac{d^i}{dx^i}\
 \frac{ A_i(\bs m,\bs l,k)\,x\,-\,B_i(\bs m,\bs l,k)}
{x^{r+1-i}\,(x-1)}\ \frac{1}{\tau(x)}\ .  
\eea
The operator has order $r+1$ and rational coefficients. Being written in the form
$$
D_{\bs m,\bs l,k}\ =\ \sum_{i=0}^{r+1} \,c_i(x)\, \frac {d^i}{dx^i}\ ,
$$
 the operator has the leading coefficient $c_{r+1}$ equal to one.

For example, for $r=1$, the operator is the classical hypergeometric differential
operator
\bea
\frac{d^2}{dx^2}\
 -\ \frac{kx + m_1(x-1)}{x(x-1)}\,\di\ +\ \frac{l_1(k+m_1+1-l_1)}{x(x-1)} \ .
\eea

\begin{theorem}
\label{main thm}
 Let $(\bs m,\bs l,k)$ be consistent parameters. Let
 $V(\bs m,\bs l,k)$ be the space of the first type associated to 
 $(\bs m,\bs l,k)$. Then the kernel of $D_{\bs m,\bs l,k}$ is $V(\bs m,\bs l,k)$.
\end{theorem}

\begin{proof}
The operator  $D_{\bs m,\bs l,k}$ is obtained from the operator 
$D^\vee_{\bs m,\bs l,k}$ by formal conjugation followed by 
the conjugation with the operator of multiplication by the function
$\tau(x)$. 

The kernel of the operator $D^\vee_{\bs m,\bs l,k}$ is the space $U(\bs
m,\bs l,k)$. Then, by standard arguments, the kernel of the operator
formally conjugated to $D^\vee_{\bs m,\bs l,k}$ consists of the
functions
\be 
\{\,\frac{W(g_1,\dots,g_r)}{T(x)}\,, \ {}   g_1,\dots,g_r\in U(\bs m,\bs l,k)\,\}\ 
\ee 
where $T(x)$ is defined in formula \Ref{REMAR}.
 
Then the kernel of the operator $D_{\bs m,\bs l,k}$ consists of the
functions 
\bea 
&&
\{\,\frac{\tau(x)\,W(g_1,\dots,g_r)}{T(x)}\,, \ {}   g_1,\dots,g_r\in U(\bs m,\bs l,k)\,\}\ = 
\phantom{aaaaaaaaaaaaaaaaa}
\\
&&
\phantom{aaaaaaaaaaaaaaaaa}
\{\,W^\dagger_U(g_1,\dots,g_r),\ {} g_1,\dots,g_r\in U(\bs
m,\bs l,k)\,\}\  .  
\eea 
This space coincides with $V(\bs m, \bs l,k)$ by
Lemma \ref{dual}.
\end{proof}

\begin{cor}
\label{COR}
 For any integers $l_1,\dots,l_r$, \
 $l_1\geq \dots \geq l_r\geq 0$, and any parameters $\bs m, k$,
 the Jacobi-Pi\~neiro 
 polynomial $P_{\bs m,\bs l,k}(x)$ is annihilated by the operator $D_{\bs m,\bs l,k}$.
\end{cor}

\begin{proof}
If the parameters $(\bs m,\bs l,k)$ 
are consistent, then the Jacobi-Pi\~neiro polynomial $P_{\bs
m,\bs l,k}(x)$ belongs to the space $V(\bs m, \bs l, k)$ by Lemma
\ref{P in V}. Therefore, the Jacobi-Pi\~neiro polynomial is
annihilated by the operator $D_{\bs m,\bs l,k}$ for consistent parameters
$(\bs m,\bs l,k)$.  This implies
Corollary \ref{COR}, since the polynomial $P_{\bs m,\bs l,k}(x)$ and
the operator $D_{\bs m,\bs l,k}$ depend on the parameters $(\bs m,\bs
l,k)$ as rational functions.
\end{proof}

\medskip

For $r=1$, the statement of the corollary is classical. It says that
the Jacobi polynomial is a solution of the hypergeometric differential
equation.  For $r=2$, see both the theorem and the corollary in
\cite{MV2}. In \cite{ABV}, a recurrent procedure in $r$ is given to
construct a differential operator annihilating the 
Jacobi-Pi\~neiro polynomial $P_{\bs m,\bs l,k}(x)$ and for $r=2$ the
operator is given explicitly.  One of the referees of this paper
pointed to us the paper \cite{CV}, in which another (or maybe the
same) differential operator annnihilating the Jacobi-Pi\~neiro
polynomial is constructed.  It would be interesting to check if the
differential operators of \cite{ABV} and \cite{CV} coincide with
$D_{\bs m,\bs l,k}$\ .

\medskip

\noindent
{\bf Remark.}
The differential operator $D_{\bs m,\bs l,k}$, 
annihilating the Jacobi-Pi\~neiro polynomial $P(\bs m,\bs l,k)$,
is uniquely determined  by the following properties:
\begin{itemize}
\item 
Coefficients of $D_{\bs m,\bs l,k}$ are rational functions of 
$x, \bs m, \bs l, k$.

\item 
If $(\bs m,\bs l,k)$ are consistent, then the kernel of
$D_{\bs m,\bs l,k}$ consists of polynomials only.

\item 
The singular points of $D_{\bs m,\bs l,k}$ are at
$x=0,1,\infty$. All singular points are regular.  The exponents of
$D_{\bs m,\bs l,k}$ at $x=0$ are $0, m_1+1, m_1+m_2+2,\dots,
m_1+\dots+m_r + r$. The exponents at $x=1$ are $0, k+1,k+2,
\dots,k+r$.  The exponents at $x=\infty$ are the numbers
$k+\sum_{s=1}^i m_s-l_{i}+l_{i+1}+i$ for $i=0,\dots,r$.
\end{itemize}
However, the
characteristic equations for exponents of $D_{\bs m,\bs l,k}$
at $x = 0, 1, \infty$ do not determine
the coefficients of $D_{\bs m,\bs l,k}$. The triviality of the
monodromy of $D_{\bs m,\bs l,k}$ 
is essential for the uniqueness of the operator
$D_{\bs m,\bs l,k}$, in
contrast with the situation for the operator
$D^\vee_{\bs m,\bs l,k}$, where the
uniqueness of  the operator $D^\vee_{\bs m,\bs l,k}$
is determined by the characteristic equations for exponents of $D^\vee_{\bs m,\bs l,k}$
only.



\begin{thebibliography}{00000}
\normalsize

\bibitem [ABV]{ABV} A. I. Aptekarev, A. Branquinho, and W. Van Assche 
{\it Multiple orthogonal polynomials for classical weights}, 
Trans. Amer. Math. Soc. 355 (2003), no. 10, 3887--3914

\bibitem[CV]{CV} J. Coussement and W. Van Assche, 
{\it Differential equations for multiple orthogonal polynomials with 
respect to classical weights: raising and lowering operators}, 
J. Phys. A: Math. Gen. 39 (2006), 3311-3318

\bibitem[MTV1]{MTV1} 
E. Mukhin, V. Tarasov, and A. Varchenko, {\it
The B. and M. Shapiro conjecture in real algebraic geometry and
the Bethe ansatz}, math.AG/0512299, 1--17

\bibitem[MTV2]{MTV2} 
E. Mukhin, V. Tarasov, and A. Varchenko, {\it
Bethe Eigenvectors of Higher Transfer Matrices}, math.QA/0605015, 1--48

\bibitem[MV1]{MV1} E. Mukhin and A. Varchenko,
{\it Critical Points of Master Functions and Flag Varieties},
Communications in Contemporary Mathematics (2004), vol. 6, no. 1, 111-163

\bibitem[MV2]{MV2} E. Mukhin and A. Varchenko, {\it Multiple orthogonal 
polynomials and a counterexample to Gaudin Bethe Ansatz Conjecture}, 
math.QA/0501144, 
1--40. To appear in Transactions of AMS

\bibitem[MV3]{MV3} E. Mukhin and A. Varchenko,
{\it Spaces of quasi-polynomials and the Bethe Ansatz},
 math.QA/0604048, 1--29




\bibitem[P]{P} L. R. Pi\~neiro, {\it On simultaneous Pade approximants
    for a collection of Markov functions}, Vestnik Mosk. Univ. Ser.,
  {\bf I}, no. 2 (1987), 52--55 (in Russian); translated in Moscow
  Univ. Math. Bull. {\bf 42}, no. 2 (1987), 52--55
  
\bibitem[ScV]{ScV} I. Scherbak and A. Varchenko, {\it Critical points
    of functions, $sl_2$ representations and Fuchsian differential
    equations with only univalued solutions}, 
Dedicated to Vladimir I. Arnold on the occasion 
of his 65th birthday.  Mosc. Math. J.  3  (2003),  no. 2, 621--645, 745

\end{thebibliography}
\end{document}